\documentclass[a4paper,12pt]{article}

\usepackage{amsthm,amsmath, amsfonts, amstext, amssymb, amsbsy}
\usepackage[dvips,pdftex]{graphicx}
\usepackage{fullpage}

\newtheorem{theorem}{Theorem}
\newtheorem{proposition}{Proposition}

\newtheorem{definition}{Definition}
\newtheorem{lemma}{Lemma}
\newtheorem{remark}{Remark}

\newcommand{\cF}{{\cal F}}
\newcommand{\cG}{{\cal G}}

\newcommand{\cL}{{\cal L}}
\newcommand{\cM}{{\cal M}}

\newcommand{\cP}{{\cal P}}

\newcommand{\cX}{{\cal X}}
\newcommand{\cY}{{\cal Y}}

\newcommand{\bx}{{\bf x}}
\newcommand{\by}{{\bf y}}
\newcommand{\bz}{{\bf z}}
\newcommand{\bt}{{\bf t}}
\newcommand{\bs}{{\bf s}}

\newcommand{\bbC}{\mathbb{C}}

\newcommand{\bbE}{\mathbb{E}}

\newcommand{\bbN}{\mathbb{N}}

\newcommand{\bbP}{\mathbb{P}}

\newcommand{\bbR}{\mathbb{R}}

\begin{document}

\title{Regular conditional distributions of max infinitely divisible processes}

\author{Cl\'ement Dombry\footnote{Universit\'e de
Poitiers, Laboratoire de Mathématiques et Applications, UMR CNRS 7348, T\'el\'eport 2, BP 30179, F-86962 Futuroscope-Chasseneuil cedex,
France.  Email: Clement.Dombry@math.univ-poitiers.fr}\ \ and Fr\'ed\'eric Eyi-Minko\footnote{Universit\'e de
Poitiers, Laboratoire de Mathématiques et Applications, UMR CNRS 7348, T\'el\'eport 2, BP 30179, F-86962 Futuroscope-Chasseneuil cedex,
France.  Email: Frederic.Eyi.minko@math.univ-poitiers.fr}}
\date{}
\maketitle

\begin{abstract}
This paper is devoted to the {\it prediction problem} in extreme
  value theory. Our main result is an explicit expression of the {\it
    regular conditional distribution} of a max-stable (or
  max-infinitely divisible) process $\{\eta(t)\}_{t\in T}$ given
  observations $\{\eta(t_i)=y_i,\ 1\leq i\leq k\}$.
  Our starting point is the point process representation of
  max-infinitely divisible processes by Gin\'e, Hahn and Vatan
  (1990). We carefully analyze the structure of the underlying point
  process, introduce the notions of {\it extremal function}, {\it
    sub-extremal function} and {\it hitting scenario} associated to
  the constraints and derive the associated distributions.  This
  allows us to explicit the conditional distribution as a mixture over
  all hitting scenarios compatible with the conditioning constraints.
  This formula extends a recent  result by Wang and Stoev
  (2011) dealing with the case of spectrally discrete max-stable
  random fields. This paper offers new tools and perspective
  for prediction in extreme value theory together with numerous
  potential applications.
\end{abstract}
\\

\section{Introduction}
\subsection{Motivations}
Since the pioneer works by Fisher and Typett \cite{FT28} and  Gnedenko
\cite{G43}, the univariate theory of extremes is now well established
with extensive studies on models, domains of attraction, parameter
estimations, {\it etc.} (see e.g. de Haan and Fereira \cite{dHF06} and
the references therein). The last decades have seen the quick
development of multivariate and spatial extreme value theory: the
emphasis is put on the characterization, modeling and estimation of
the dependence structure of multivariate extremes. Among many others,
the reader should refer to the excellent monographs
\cite{BGT04,dHF06,EKM97,R08} and the reference therein.

Max-stable random fields turn out to be fundamental models for spatial extremes since they arise as the the limit of rescaled
maxima. More precisely, consider the component-wise maxima
\[
\eta_n(t)=\max_{1\leq i\leq n} X_i(t),\quad t\in T,
\]
of independent and identically distributed (i.i.d.) random fields  $\{X_i(t)\}_{t\in T}$, $i\geq 1$.
 If the random field
$\eta_n= \{ \eta_n(t) \}_{t\in T}$ converges in distribution, as
$n\to\infty$, under suitable affine normalization, then its limit
$\eta=\{\eta(t)\}_{t\in T}$ is necessarily max-stable (see
e.g. \cite{dHF06,R08}). Therefore, max-stable random fields play a
central role in extreme value theory, just like Gaussian random fields
do in the classical statistical theory based on the Central Limit
Theorem.

In this framework, the {\it prediction problem} arises as an important
and long-standing challenge in extreme value theory. Suppose that we
already have a suitable max-stable model for the dependence structure
of a random field $\eta=\{\eta(t)\}_{t\in T}$ and that the field is
observed at some locations $t_1,\ldots,t_k\in T$. How can we take
benefit from these observations and predict the random field $\eta$ at
other locations ?
 We are naturally lead to consider the {\it
  conditional distribution} of $ \{ \eta(t) \}_{t\in T}$ given the
 observations $\{\eta(t_i)=y_i,\ 1\leq i\leq k\}$.
A formal definition of the notion of regular conditional distribution
is  deferred to the Appendix~\ref{sec:rcd}.

In the classical Gaussian framework, i.e., if $\eta$ is a Gaussian
random field, it is well known that the corresponding conditional
distribution remains Gaussian and simple formulas give the conditional
mean and covariance structure. This theory is strongly linked with the
theory of Hilbert spaces: the conditional expectation, for example, can be obtained
as the $L^2$-projection of the random field $\eta$ onto a
suitable Gaussian subspace.  In extreme value theory, the prediction
problem turns out to be difficult. A first approach  by Davis and Resnick \cite{DR89,DR93}
is based on a $L^1$-metric between max-stable variables and on a kind of projection onto 
max-stable spaces. To some extent, this work mimics the corresponding $L^2$-theory for Gaussian
spaces. However, unlike the Gaussian case, there is no clear relationship between the  predictor obtained by projection onto the max-stable space generated by the variables $\{\eta(t_i),\ 1\leq i\leq k\}$ and the conditional distributions of $\eta$ with respect to these variables. A first  major contribution to the conditional distribution
problem is the work by Wang and Stoev \cite{WS11}. The authors consider max-linear random fields, a special
class of max-stable random fields with discrete spectral measure, and
give an exact expression of the conditional distributions as well as
efficient algorithms. The max-linear structure plays an essential role
in their work and  provides major simplifications since in this case $\eta$ admits the simple representation
\begin{equation}\label{eq:maxlin}\nonumber
\eta(t)=\bigvee_{j=1}^q Z_j f_j(t), \quad t\in T,
\end{equation}
where the symbol $\bigvee$ denotes the maximum, $f_1,\ldots,f_q$ are deterministic functions and $Z_1,\ldots,Z_q$ are i.i.d. random variables with unit Fr\'echet distribution. The authors  determine  the conditional distributions of $(Z_j)_{1\leq j\leq q}$ given observations $\{\eta(t_i)=y_i,\ 1\leq i\leq k\}$. Their result relies on the important notion of {\it hitting scenario} defined as the subset  of indices  $j\in[\![1,q]\!]$ such that $\eta(t_i)=Z_jf(t_i)$ for some $i\in [\![1,k]\!]$, where, for $n\geq 1$, we note $[\![1,n]\!]=\{1,\ldots,n\}$. The conditional distribution of $(Z_j)_{1\leq j\leq q}$ is expressed as a mixture over all admissible hitting scenarios with minimal rank.

The purpose of the present  paper is to propose a general theoretical framework for 
conditional distributions in extreme value theory,  covering not only
the whole class of sample continuous max-stable random fields but also
the class of sample continuous max-infinitely divisible (max-i.d.)
random fields (see Balkema and Resnick \cite{BR77}).
Our starting point is the general representation by Gin\'e, Hahn and Vatan
\cite{GHV90} of  max-i.d. sample continuous random fields (see also de Haan \cite{dH84} for the max-stable case). 
It is possible to construct a Poisson random measure $\Phi=\sum_{i=1}^N \delta_{\phi_i}$ on the space of continuous functions on $T$ such that
\[
\eta(t)\stackrel{\cL}=\bigvee_{i=1}^N \phi_i(t),\quad t\in T.
\]
Here the random variable $N$ is equal to the total mass of $\Phi$ that  may be finite or infinite and  $\stackrel{\cL}=$ stands for equality of probability laws (see Theorem \ref{theo:GHV} below for a precise statement). We denote by $[\Phi]=\{\phi_i,\ 1\leq i\leq N\}$ the set of atoms of $\Phi$. Clearly, $\phi(t)\leq \eta(t)$ for all $t\in T$ and $\phi\in[\Phi]$. The observations $\{\eta(t_i)=y_i,\ 1\leq i\leq k\}$ naturally lead to consider extremal points: a function $\phi\in[\Phi]$ is called {\it extremal} if $\phi(t_i)=\eta(t_i)$ for some $i\in[\![1,k]\!]$,
otherwise it is called {\it sub-extremal}. 
We show that under some mild condition, one can  define a random partition $\Theta=(\theta_1,\ldots,\theta_\ell)$ of $\{t_1,\ldots,t_k\}$ and extremal functions $\varphi_1^+,\ldots,\varphi_\ell^+ \in[\Phi]$ such that the point $t_i$ belongs to the component $\theta_j$  if and only if $\varphi_j^+(t_i)=\eta(t_i)$. Using the terminology of Wang and Stoev \cite{WS11},
we call  {\it hitting scenario} a partition of $\{t_1,\cdots,t_k\}$ that reflects the way how  the extremal functions $\varphi_1^+,\ldots,\varphi_\ell^+$ hit the constraints $\varphi_j^+(t_i)\leq \eta(t_i)$, $1\le i\leq k$. The main results of this paper are Theorems~\ref{theo2.2} and~\ref{theo3}, where the conditional distribution of $\eta$ given $\{\eta(t_i)=y_i,\ 1\leq i\leq k\}$ is expressed as a {\it mixture} over all possible hitting scenarios.

The paper is structured as follows. In Section 2, the distribution of
extremal and sub-extremal functions is analyzed and a characterization of
the hitting scenario distribution is given.  In Section 3, we focus on conditional distributions:  we compute the conditional distribution of the hitting scenario and extremal functions  and then derive the conditional distribution of $\eta$. Section 4 is devoted to examples: we specify our results in the simple case of a single conditioning point and consider  max-stable models. The proofs are collected in Section 5 and some technical details are postponed to an appendix.

\subsection{Preliminary on max-i.d. processes}
Let $T$ be a compact metric space and $\bbC=\bbC(T,\bbR)$ be the space of continuous functions on $T$ endowed with the sup norm
\[
 \|f\|=\sup_{t\in T} |f(t)|,\quad f\in\bbC.
\]
Let $(\Omega,\cF,\bbP)$ be a  probability space. A random
process $\eta=\{\eta(t)\}_{t\in T}$ is said to be max-i.d. on $\bbC$
if $\eta$ has a version with continuous sample path and if, for each
$n\geq 1$, there exist $\{\eta_{ni},\ 1\leq i\leq n\}$ i.i.d. sample continuous random fields on
$T$ such that
\[
 \eta \stackrel{\cL}= \bigvee_{i=1}^n \eta_{ni},
\]
where $\bigvee$ denotes pointwise maximum.

Gin\'e, Hahn and Vatan (see \cite{GHV90} Theorem 2.4) give a
representation of such processes in terms of Poisson random
measure. For any function $f$ on $T$ and set $A\subset T$, we note
$f(A)=\sup_{t\in A} f(t)$.
\begin{theorem}(Gin\'e, Hahn and Vatan \cite{GHV90})\label{theo:GHV}\\
Let $h$ be the {\it vertex function} of a sample continuous max-i.d. process $\eta$ defined by
\[
 h(t)=\sup\{x\in\bbR;\ \bbP(\eta(t)\geq x)=1\}\in [-\infty,\infty),\quad t\in T,
\]
and define $\bbC_h=\{f\in\bbC;\ f\neq h, f\geq h\}$. Under the condition that the vertex function $h$ is continuous, there exists a locally-finite Borel measure $\mu$ on $\bbC_h$, such that if $\Phi$ is  a Poisson random measure $\Phi$ on $\bbC_h$ with intensity measure $\mu$, then
\begin{equation}
\left\{ \eta(t)\right\}_{t\in T} \stackrel{\cL} = \left\{\sup\{h(t),\phi(t); \phi\in [\Phi]\}\right\}_{t\in T} \label{repGHV}
\end{equation}
where $[\Phi]$ denotes the set of atoms of $\Phi$.\\ 
Furthermore, the following relations hold:
\begin{equation}
  \label{eq:C1}
 h(K)=\sup\{x\in\bbR;\ \bbP(\eta(K)\geq x)=1\},\quad K\subset T\ \mathrm{closed}, 
\end{equation}
and
\begin{equation}
  \label{eq:C2}
  \bbP\left[\eta(K_i)<x_i,\ 1\leq i\leq  n \right]=\exp  [
    -\mu \left( \cup_{i=1}^n \{f\in\bbC_h;\ f(K_i)\geq x_i\} \right)],
\end{equation}
where $n\in\bbN$, $K_i\subset T$  closed and $x_i>h(K_i)$, $1\leq i\leq n$.
\end{theorem}

Theorem~\ref{theo:GHV} provides an almost complete description of
max-i.d. continuous random processes, the only restriction being the
continuity of the vertex function.  Clearly, the distribution of
$\eta$ is completely characterized by the vertex function $h$ and the
so called {\it exponent measure} $\mu$. The random process  
$e^\eta-e^h$ is continuous and max-i.d. and its vertex function is identically equal to $0$.
%
%
%
%
Since the conditional distribution of $\eta$ is easily deduced from
that of $e^\eta-e^h$, we can assume without loss of generality that $h \equiv
0$; the corresponding set $\bbC_0$ is the space of non negative and
non null continuous functions on $T$.

We need some more notations from point process theory (see Daley and
Vere-Jones \cite{DVJ03,DVJ08}).  It will be convenient to introduce  
a measurable enumeration of the atoms of $\Phi$ (see \cite{DVJ08} Lemma 9.1.XIII).
 The total mass of $\Phi$ is noted $N=\Phi(\bbC_0)$. If  $\mu(\bbC_0)<\infty$, $N$ has  a Poisson distribution with mean $\mu(\bbC_0)$,
otherwise $N=+\infty$  almost surely (a.s.). One can construct $\bbC_0$-valued 
random variables $(\phi_i)_{i\geq 1}$ such that $\Phi=\sum_{i=1}^{N}\delta_{\phi_i}$.

Let $M_p(\bbC_0)$ be the space of point measures $M=\sum_{i\in I}\delta_{f_i}$ on $\bbC_0$ such that 
\[
\{f_i \in \bbC_0\colon \|f_i\| > \varepsilon \}\ \mathrm{is\ finite\ for\ all}\  \varepsilon > 0.
\]
We endow $M_p(\bbC_0)$ with the $\sigma$-algebra $\cM_p$
generated by the applications 
\[
 M\mapsto M(A),\quad A\subset \bbC_0\mbox{\ Borel set\ }.
\]
For $M\in M_p(\bbC_0)$, let $[M]=\{f_i,\ i\in I\}$ be the countable set of 
atoms of $M$.  If $M$ is non null, then for all $t\in T$, the set
$\{f(t); f\in [M]\}$ is non empty and has finitely many points in $(\varepsilon,+\infty)$ for all $\varepsilon>0$ so that the maximum $\max\{f(t); f\in [M]\}$ is reached. Furthermore by considering restrictions of the measure $M$ to
sets $\{f\in\bbC_0; \ \|f\|>\varepsilon\}$ and using  uniform
convergence, it is easy to show that the mapping
\[
\max(M)\colon
\begin{cases}
  T \rightarrow [0,+\infty) \\
  t \mapsto \max\{f(t); f\in [M]\}
\end{cases}
\]
is continuous with the convention that $\max(M)\equiv 0$ if $M=0$. 

In Theorem \ref{theo:GHV} (with $h\equiv 0$), Equation \eqref{eq:C2} implies that the exponent measure
$\mu$ satisfies, for all $\varepsilon>0$,
\begin{equation}\label{eq:condmu}
\mu(\{f\in\bbC_0;\ \|f\|>\varepsilon\})<\infty.  
\end{equation}
Consequently, we have $\Phi\in M_p(\bbC_0)$ almost surely and $\eta\stackrel{\cL}=\max(\Phi)$. 

\begin{figure}[t]  
\begin{center}
\includegraphics[width=5cm]{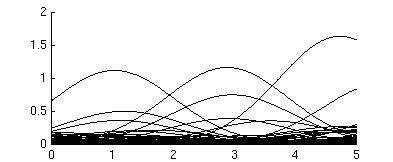}
  \includegraphics[width=5cm]{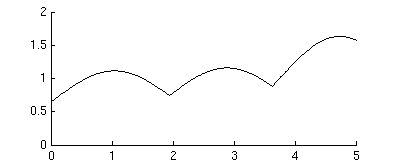}
\end{center}    
\caption{A representation of the point process $\Phi$ (left) and of the associated maximum process $\eta=\max(\Phi)$ (right) in the moving maximum max-stable model based on the gaussian density function. Here $T=[0,5]$.}
\label{fig1}
\end{figure}

An illustration of Theorem \ref{theo:GHV} is given in Figure \ref{fig1} with a representation of the Poisson point measure $\Phi$ and  of the corresponding maximum process $\eta=\max(\Phi)$ in the moving maximum max-stable model based on the Gaussian density function. 
 
\section{Extremal points and related distributions}

In the sequel, $\eta$ denotes a sample continuous max-i.d. random
process with vertex function $h\equiv 0$ and exponent measure $\mu$ on
$\bbC_0$. On the same probability space, we suppose that a
$M_p(\bbC_0)$-valued Poisson random measure $\Phi=\sum_{i=1}^N\delta_{\phi_i}$ with intensity
measure $\mu$ is given and such that $\eta=\max(\Phi)$.

\subsection{Extremal and sub-extremal point measures}

Let $K\subset T$ be a closed subset of $T$. We introduce here the
notion of $K$-extremal points that will play a key role in this
work. We use the following notations: if $f_1,f_2$ are two functions
defined (at least) on $K$, we write

\begin{align*}
  f_1=_Kf_2\quad&\text{if and only if}\quad \forall t\in K,\
  f_1(t)=f_2(t),\\
  f_1<_Kf_2\quad&\text{if and only if}\quad \forall t\in K,\
  f_1(t)<f_2(t),\\
  f_1\not<_Kf_2\quad&\text{if and only if}\quad \exists t\in K,\
  f_1(t)\geq f_2(t).
\end{align*}

Let $M\in M_p(\bbC_0)$. An atom $f\in [M]$ is called $K$-sub-extremal if and only if $f<_K \max(M)$ and $K$-extremal otherwise. In words, a sub-extremal atom has no contribution to the maximum $\max(M)$ on $K$.

\begin{definition}\label{def:ep}
Define the $K$-extremal random point measure $\Phi_K^+$ and the $K$-sub-extremal random point measure $\Phi_K^-$ by
\[
\Phi_K^+=\sum_{i=1}^N 1_{\{\phi_i\not<_K \eta\}}\delta_{\phi_i}\quad \mbox{and}\quad \Phi_K^-=\sum_{i=1}^N 1_{\{\phi_i<_K \eta\}}\delta_{\phi_i}.
\]
\end{definition}

\begin{figure*}[t]  
\begin{center}  
  \includegraphics[width=5cm]{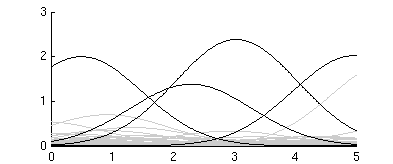}
  \includegraphics[width=5cm]{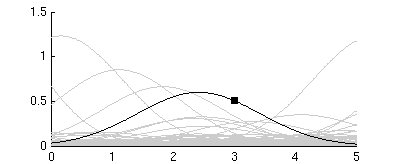}
\end{center}
  \caption{Decomposition of the Poisson point measure $\Phi$ into the $K$-extremal point measure $\Phi_K^+$ (black) and the $K$-sub-extremal point measure $\Phi_K^-$ (grey). Left:  $K=[0,5]$. Right:  $K=\{3\}$ represented by a black square.}
\label{fig2}
\end{figure*}

Figure \ref{fig2} provides an illustration of the definition. It should be noted that   $\Phi_K^+$ and $\Phi_K^-$ are well defined measurable random point measures (see Lemma~\ref{lem:meas1} in  Appendix~\ref{sec:meas}). 
Furthermore, it is straightforward from the definition that
\[
\Phi= \Phi_K^+ + \Phi_K^-,\quad \max(\Phi_K^+)=_K \eta \quad  \mbox{and} \quad \max(\Phi_K^-)<_K \eta.  
\]
Define the following measurable  subsets of $M_p(\bbC_0)$ (see Lemma~\ref{lem:meas2} in Appendix~\ref{sec:meas}):
\begin{eqnarray}
C^+_K&=&\Big\{M\in M_p(\bbC_0);\ \forall f\in [M],\ f \not<_K  \max(M)\Big\}\label{eq:Csharp},\\
C^-_K(g)&=&\Big\{M\in M_p(\bbC_0);\ \forall f\in [M],\   f<_K g\Big\}\label{eq:Cflat},
\end{eqnarray}
where $g$ is any continuous function defined (at least) on $K$.
Clearly, it always holds 
\[
\Phi_K^+ \in C^+_K \quad \mathrm{and}\quad  \Phi_K^-\in C^-_K(\eta).
\]

The following theorem characterizes the joint distribution of $(\Phi_K^+,\Phi_K^-)$ given that $\Phi_K^+$ is
finite. We note $\delta_0$ the Dirac mass at $0$.
\begin{theorem}\label{theo1}
For all measurable $A,B\subset M_p(\bbC_0)$,
\[
\bbP\big[(\Phi_K^+,\Phi_K^-)\in A\times B  ,\ \Phi_K^+(\bbC_0)=0\big]
=\exp[-\mu(\bbC_0)]\,\delta_0(A)\,\delta_0(B),
\]
and, for $k\geq 1$,
\begin{eqnarray*} 
&&\bbP\big[(\Phi_K^+,\Phi_K^-)\in A\times B  ,\ \Phi_K^+(\bbC_0)=k\big]\\
&=&\frac{1}{k!} \int_{\bbC_0^k}1_{\{\sum_{i=1}^k \delta_{f_i}\in A\cap C^+_K\}} \bbP\left[ \Phi \in B \cap C^-_K \left(  \vee ^k_{i=1} f_i \right)  \right] \,\mu^{\otimes k}(df_1,\ldots,df_k).
\end{eqnarray*}
\end{theorem}

Theorem \ref{theo1} fully  characterizes the joint distribution of $(\Phi_K^+,\Phi_K^-)$ provided that $\Phi_K^+(\bbC_0)$ is almost surely finite. We now focus on this last condition.

\begin{proposition}\label{prop1}\ \\
The $K$-extremal point measure $\Phi_K^+$ is a.s. finite if and only if one of the following condition holds:
\begin{itemize}
\item[{\rm (i)}]  $\mu(\bbC_0)<+\infty$;
\item[{\rm (ii)}] $\mu(\bbC_0)=+\infty$  {\rm and}  $\inf\limits_{t\in K} \eta(t)>0$  {\rm almost surely}.
\end{itemize}
\end{proposition}

It should be noted that any simple max-stable random field (with unit Fr\'echet margins) satisfies condition (ii) above. See for example Corollary 3.4 in \cite{GHV90}. 

\begin{remark}
Using Theorem \ref{theo1}, it is easy to show that the distribution of $(\Phi_K^+,\Phi_K^-)$ has the following structure. Define  the tail functional $\bar\mu_K$  by
\[
\bar\mu_K(g)=\mu(\{f\in\bbC_0;\ f\not<_K g\})
\]
for any continuous function $g$  defined (at least) on  $K$. 
Suppose that $\Phi_K^+$ is  finite almost surely. Its distribution is then given by the so-called Janossy measures (see e.g. Daley and Vere-Jones \cite{DVJ03} section 5.3). The Janossy measure of order $k$ of the $K$-extremal point measure $\Phi_K^+$ is given by
\[
J_k(df_1,\ldots,df_k)=\exp\left[-\bar\mu_K \left( \vee_{i=1}^k f_i \right)  \right] \, 1_{ \left\lbrace \sum_{i=1}^k \delta_{f_i}\in  C^+_K \right\rbrace }\mu^{\otimes k}(df_1,\ldots,df_k).
\]
Furthermore, given that $\Phi_K^+=\sum_{i=1}^k \delta_{f_i}$, the conditional distribution of $\Phi_K^-$ is equal to the distribution of a Poisson random measure with measure intensity  $1_{ \{f <_K \vee_{i=1}^k f_i \}}\mu(df)$. These results are not used in the sequel and we omit their proof for the sake of brevity.

\end{remark}

\subsection{Extremal functions}

Let $t\in T$. We denote by $\mu_{t}$  the measure on $(0,+\infty)$ defined by
\[
\mu_{t}(A)=\mu(\{f\in\bbC_0; f(t)\in A\}),\quad A\subset (0,+\infty)\ \mathrm{Borel\ set},
\]
and by $\bar\mu_t$ the associated tail function defined by
\[
\bar\mu_t(x)=\mu_t([ x,+\infty)),\quad x>0.
\]
Note that
\begin{equation}\label{eq:C2bis}
\bbP(\eta(t) < x) = \exp(-\mu(\{f\in\bbC_0; f(t) \geq x\}))=\exp(-\bar\mu_t(x)),\quad x >0.
\end{equation}
The following proposition states that, under a natural condition, 
there is almost surely a unique $\{t\}$-extremal point in $\Phi$. This extremal
point will be referred to as the ${t}$-{\it extremal function} and noted $\phi_t^+$.
\begin{proposition}\label{prop2}
For $t\in T$, the following statements are equivalent:
\begin{enumerate}
\item[{\rm (i)}] $\Phi_{\{t\}}^+(\bbC_0)=1$ almost surely;
\item[{\rm (ii)}] $\bar\mu_{t}(0^{+})=+\infty$ and $\bar\mu_t$ is continuous on $(0,+\infty)$;
\item[{\rm (iii)}] the distribution of $\eta(t)$ has no atom.
\end{enumerate}
If these conditions are met, we define  the $t$-extremal function $\phi_t^+$
 by the relation  $\Phi_{\{t\}}^+=\delta_{\phi^+_t}$ a.s..  For all
measurable $A\subset \bbC_0$ we have
\begin{eqnarray}
\bbP(\phi_t^+\in A)&=&\int_A \exp[-\bar\mu_{t}(f(t))]\,\mu(df).  \label{eq:prop2.2}
\end{eqnarray}
\end{proposition}
An important class of processes satisfying the conditions of Proposition~\ref{prop2} is the class of max-stable processes (see section \ref{sec:maxstable} below).
 
\subsection{Hitting scenarios} 
Proposition~\ref{prop2} gives the distribution of $\Phi_K^+$ when $K=\{t\}$ is reduced to a single point. 
Going a step further, we consider the case when $K$ is finite. In the sequel, we suppose that the following assumption is satisfied:
\begin{center}
 {\rm({\bf A})}\ \  $K=\{t_1,\ldots,t_k\}$ is finite and, \\for all $t\in K$, $\bar\mu_{t}$ is continuous and $\bar\mu_{t}(0^{+})=+\infty$.
\end{center}
Roughly speaking, this  ensures that the maximum $\eta(t)=\max(\Phi)(t)$ is uniquely reached for all $t\in K$. This will provide combinatorial simplifications. More precisely, under Assumption ${\rm({\bf A})}$, the event 
$$\Omega_K=\bigcap_{t\in K}\{\Phi_{\{t\}}^{+}(\bbC_0)=1\}$$ is of probability $1$ and  the extremal functions $\phi_{t_1}^+,\ldots,\phi_{t_k}^+$ are  well defined. In the next definition, we introduce the notion of hitting scenario that reflects the way how these extremal functions  hit the maximum $\eta$ on $K$.

Let $\cP_K$ be the set of partitions of $K$. It is  convenient to think about $K$ as an ordered set, say $t_1<\cdots<t_k$.  Then each partition $\tau$ can be written uniquely in the standardized form $\tau=(\tau_1,\ldots,\tau_\ell)$ where $\ell=\ell(\tau)$ is the length of the partition, $\tau_1\subset K$ is the component of $t_1$, $\tau_2\subset K$ is the component containing $\min(  K\setminus\tau_1)$ and so on. With this convention, the components $\tau_1,\ldots,\tau_\ell$ of the partition are labeled so that $$\min \tau_1<\cdots <\min\tau_\ell.$$

\begin{figure}[t]  
\begin{center}  
\includegraphics[width=5cm]{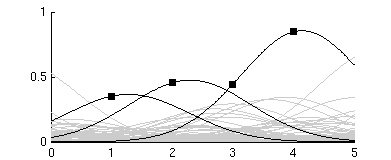}
 \includegraphics[width=5cm]{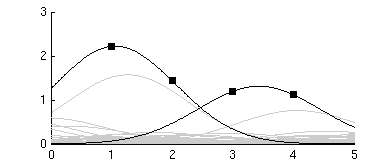}
\end{center}
\caption{Two realisations of the Poisson point measure $\Phi$ and of the corresponding hitting scenario $\Theta$ and extremal functions $\varphi_1^+,\ldots,\varphi_{\ell(\Theta)}^+$ with  $K=\{t_1,t_2,t_3,t_4\}$  represented by the black squares. Left: the hitting scenario is $\Theta=(\{t_1\},\{t_2\},\{t_3,t_4\})$, the extremal functions  are $\varphi_1^+=\phi_{t_1}^+$, $\varphi_2^+=\phi_{t_2}^+$ and $\varphi_3^+=\phi_{t_3}^+=\phi_{t_4}^+$. Right: the hitting scenario is $\Theta=(\{t_1,t_2\},\{t_3,t_4\})$, the extremal functions (black) are $\varphi_1^+=\phi_{t_1}^+=\phi_{t_2}^+$ and $\varphi_2^+=\phi_{t_3}^+=\phi_{t_4}^+$.
}
\label{fig3}
\end{figure}

\begin{definition}\label{def:hs}
Suppose that Assumption ${\rm({\bf A})}$ is met.
Define $\sim$ the (random) equivalence relation on $K=\{t_1,\ldots,t_k\}$  by $t\sim t'\quad \mbox{if\ and\ only\ if}\quad \phi_{t}^+=\phi_{t'}^+$.
The  partition $\Theta=(\theta_1,\ldots,\theta_{\ell(\Theta)})$ of $K$ into equivalence classes is called the hitting scenario.
For $j\in[\![1,\ell(\Theta)]\!]$, let $\varphi_j^+$ be the extremal function associated to the component $\theta_j$, i.e., such that $\varphi_{j}^+=\phi^+_{t}$ for all $t\in \theta_j$.  
\end{definition}
We illustrate the definition with two examples in Figure~\ref{fig3}.
Clearly a point $\phi\in[\Phi]$ is $K$-extremal if and only if it is ${t}$-extremal for some $t\in K$, so that  $[\Phi_K^+]=\{\phi_t^+,\ t\in K\}$. Furthermore, the random measure $\Phi_K^+$ is almost surely simple, i.e. any atoms have a simple multiplicity, otherwise the condition $\Phi_{\{t\}}^{+}(\bbC_0)=1$ a.s. would not be satisfied for some $t\in K$. These considerations entail that 
\begin{equation}\label{eq:decPhiK+}
\Phi_K^+=\sum_{j=1}^{\ell(\Theta)} \delta_{\varphi_j^+}.
\end{equation}
In particular, the length $\ell(\Theta)$ of the hitting scenario is equal to $\Phi_K^+(\bbC_0)$. Furthermore the extremal functions satisfy
\begin{equation}
 \forall j\in [\![1,\ell]\!],\ \forall t\in\theta_j,\quad \varphi_{j}^+(t)  >  \bigvee_{j'\neq j} \varphi_{j'}^+(t). \label{eq:ineqhs}
\end{equation}

The distribution of the hitting scenario and extremal functions is given by the following
proposition. The proof relies on Theorem \ref{theo1}.

\begin{proposition}\label{prop3} Suppose Assumption {\rm({\bf A})} is met.\\
  Then, for any partition $\tau=(\tau_1,\ldots,\tau_\ell)\in\cP_K$, and any Borel sets
  $A\subset\bbC_0^{\ell}$, $B\subset M_p(\bbC_0)$, we have
 \begin{eqnarray}
&&\bbP[\Theta=\tau,\ (\varphi^+_{1},\ldots,\varphi^+_{\ell})\in A,\ \Phi_K^-\in B]\label{eq:prop3.1}\\
&=& \int_{A}1_{\{\forall j\in[\![1,\ell]\!],\ f_j >_{\tau_j} \bigvee_{j'\neq j} f_{j'} \}} 
\bbP \left[ \Phi \in B \cap C^-_K \left( \vee^\ell_{j=1} f_j \right) \right]
\,\mu^{\otimes \ell}(df_1,\ldots,df_\ell). \nonumber
\end{eqnarray}
\end{proposition}

\section{Regular conditional distribution of max-id processes}

We now focus on  conditional distributions. We will need some notations. 

If $\bs=(s_1,\ldots,s_l)\in T^l$ and $f\in \bbC_0$, we note
$f(\bs)=(f(s_1),\ldots,f(s_l))$.  Let $\mu_{\bs}$ be the exponent measure  of the max-i.d. random vector $\eta(\bs)$, i.e. the measure on $[0,+\infty)^l\setminus\{0\}$ defined by
\[
 \mu_\bs(A)=\mu(\{f\in\bbC_0;\ f(\bs)\in A\}),\quad A\subset [0,+\infty)^l\setminus\{0\} \mbox{\ Borel\ set.}
\]
Define the corresponding tail function
\[
\bar\mu_\bs(\bx)=\mu(\{f\in\bbC_0;\ f(\bs)\not < \bx\}),\quad \bx\in [0,+\infty)^l \smallsetminus \{ 0 \}.
\]  
Let $\{P_{\bs}(\bx,df);\ \bx\in[0,+\infty)^l\setminus\{0\}\}$ be a regular version of the conditional measure $\mu(df)$ given $f(\bs)=\bx$ (see Lemma~\ref{lem:rcd} in Appendix~\ref{sec:rcd}). Then for any measurable function $$F: [0,+\infty)^{l}\times\bbC_0\to\ [0,+\infty)$$ vanishing on $\{0\}\times\bbC_0$, we have
\begin{eqnarray}
\int_{\bbC_0} F(f(\bs),f)\mu(df) & =&\int_{[0,+\infty)^l\setminus\{0\}}\int_{\bbC_0} F(\bx,f)P_{\bs}(\bx,df)\mu_{\bs}(d\bx).\label{eq:theo2.1pre}
\end{eqnarray}

Let $\bt=(t_1,\ldots,t_k)$ and $\by=(y_1,\ldots,y_k)\in [0,+\infty)^k$. Before considering the conditional distribution of $\eta$ with respect to $\eta(\bt)=\by$, we give in the next theorem  an explicit expression of the distribution  of $\eta(\bt)$. We note $ K=\{t_1,\ldots,t_k\}$. For any non empty $L\subset K$, we define  $\tilde L= \{ i \in [\![1,k]\!] : t_i\in L\}$ and  set $\bt_L=(t_i)_{i\in \tilde L}$, $\by_L=(y_i)_{i\in \tilde L}$ and $L^{c} = K\smallsetminus L$. 
\begin{theorem}\label{theo2.1}
Suppose assumption {\rm({\bf A})} is satisfied. For $\tau\in\cP_K$, define the measure $\nu_{\bt}^\tau$ on $[0,+\infty)^k$ by
\[
\nu_{\bt}^\tau(C)=\bbP(\eta(\bt)\in C;\ \Theta=\tau),\quad C\subset [0,+\infty)^k\ \mathrm{Borel\ set}.
\]
Then,  
\begin{eqnarray}
\nu_{\bt}^\tau(d\by)
&=&  \exp[-\bar\mu_{\bt}(\by)]\bigotimes_{j=1}^\ell \Big\{P_{\bt_{\tau_j}}(\by_{\tau_j},\{ f(\bt_{\tau_j^c})< \by_{\tau_j^c}\} )\  \mu_{\bt_{\tau_j}}(d\by_{\tau_j})\Big\}\label{eq:theo2.1}
\end{eqnarray}
and the distribution $\nu_{\bt}$ of $\eta(\bt)$ is equal to $\nu_\bt=\sum_{\tau\in\cP_K}\nu_{\bt}^\tau$. 
\end{theorem}
Under some extra regularity assumptions, one can even get an explicit density function for $\nu_\bt$ (see the section~\ref{sec:rm} on regular models below).

We are now ready to state our main result. In Theorem \ref{theo2.2} below, we consider the regular conditional distribution of the point process $\Phi$ with respect to $\eta(\bt)=\by$. Then, thanks to the relation $\eta=\max(\Phi)$, we deduce easily in Corollary \ref{theo3} below the regular conditional distribution of $\eta$ with respect to $\eta(\bt)=\by$. 

Recall that the point process has been decomposed into two parts: a hitting scenario $\Theta$ together with extremal functions $(\varphi_1^+,\ldots,\varphi_{\ell(\Theta)}^+)$ and a $K$-sub-extremal point process $\Phi_K^-$. Taking this decomposition into account, we introduce the following regular conditional distributions:
\begin{eqnarray*}
\pi_\bt(\by,\,\cdot\,)&=& \bbP[\Theta\in\cdot \mid \eta(\bt)=\by] \\
Q_\bt(\by,\tau,\,\cdot\,)&=& \bbP[(\varphi_j^+)\in \cdot \mid \eta(\bt)=\by, \Theta=\tau] \\
R_\bt(\by,\tau,(f_j),\,\cdot\,)&=& \bbP[\Phi_K^-\in \cdot \mid \eta(\bt)=\by, \Theta=\tau, (\varphi_j^+)=(f_j)]. 
\end{eqnarray*}
We use here the short notations $\ell=\ell(\tau)$, $(\varphi_j^+)=(\varphi_1^+,\ldots,\varphi_\ell^+)$ and similarly $(f_j)=(f_1,\ldots,f_\ell)$. 
The following theorem provides explicit expressions for these regular conditional distributions.

\begin{theorem}\label{theo2.2}
Suppose assumption {\rm({\bf A})} is satisfied. 
\begin{enumerate}
\item  For any $\tau\in\cP_K$, it holds $\nu_\bt(d\by)$-a.e.
\begin{equation}\label{eq:theo2.2}
\pi_\bt(\by,\tau)=\frac{d\nu_\bt^\tau}{d\nu_\bt}(\by)
\end{equation}
where $\nu_\bt$ and $\nu_\bt^\tau$ are defined in Theorem \ref{theo2.1} and $d\nu_\bt^\tau / d\nu_\bt$ denotes the Radon-Nykodym derivative of $\nu_\bt^\tau$ w.r.t. $\nu_\bt$.
\item It holds $\nu_\bt(d\by)\pi_\bt(\by,d\tau)$-a.e. 
\begin{equation}
Q_\bt(\by,\tau,df_1\cdots df_\ell)=  \bigotimes_{j=1}^\ell \Big\{
  \frac { 1_{\{  f_j(\bt_{\tau_j^c})< \by_{\tau_j^c}\}}P_{\bt_{\tau_j}}(\by_{\tau_j},df_j)} { P_{\bt_{\tau_j}}(\by_{\tau_j},\{  f(\bt_{\tau_j^c})< \by_{\tau_j^c}\} )}   \Big\}.\label{eq:theo2.3}
\end{equation}
In words, conditionally on $\eta(\bt)=\by$ and $\Theta=\tau$, the extremal functions $(\varphi_1^+,\ldots,\varphi_\ell^+)$ are independent and $\varphi_j^+$ follows the distribution $P_{\bt_{\tau_j}}(\by_{\tau_j},df)$ conditioned to the constraint $f(\bt_{\tau_j^c})< \by_{\tau_j^c}$.
\item Let $C_\bt^-(\by)=\{M\in M_p(\bbC_0);\ \forall f\in [M],\ f(\bt)<\by\}$. It holds a.e.
\begin{equation}
R_\bt(\by,\tau,(f_j),B)\equiv R_\bt(\by,B)= \frac{\bbP[\Phi \in B \cap C_\bt^-(\by)]}{\bbP[\Phi \in C_\bt^-(\by)]}\label{eq:theo2.4}
\end{equation}
for any measurable $B\in M_p(\bbC_0)$. In words, conditionally on $\eta(\bt)=\by$, $\Phi_K^-$ is independent of $\Theta$ and $(\varphi_1^+,\ldots,\varphi_{\ell(\Theta)}^+)$ and has the same distribution as a Poisson point measure with intensity $1_{\{f(\bt)<\by\}}\mu(df)$.
\end{enumerate}
\end{theorem}

As a consequence, we deduce the regular conditional distribution of $\eta$ with respect to $\eta(\bt)=\by$.
\begin{theorem}\label{theo3}
It holds $\nu_\bt(d\by)$-a.e.
\begin{eqnarray*}
&&\bbP[\eta(\bs) < \bz  \mid \eta(\bt)=\by]= \exp[-\mu(\{ f(\bs)\not<\bz, \  f(\bt)< \by\})] \\
&&\qquad\times \sum_{\tau\in\cP_K}\pi_\bt(\by,\tau) \prod_{j=1}^\ell \frac{P_{\bt_{\tau_j}}( \by_{\tau_j},\{ f(\bt_{\tau_j^c})< \by_{\tau_j^c}, \ \ f(\bs)<\bz\})}{P_{\bt_{\tau_j}}(\by_{\tau_j},\{ f(\bt_{\tau_j^c})< \by_{\tau_j^c}\} )} 
\end{eqnarray*}
for any $l\geq 1$, $\bs\in T^l$ and $\bz\in [0,+\infty)^l$.
\end{theorem}

\begin{remark}{\rm
Let us mention that Theorem \ref{theo2.2} suggests a three-step procedure for sampling from the conditional distribution of $\eta$ given $\eta(\bt)=\by$:
\begin{enumerate}
 \item Draw a random partition $\tau$ with distribution $\pi_{\bt}(\by,\cdot)$.
 \item Given $\tau=\{\tau_1,\ldots,\tau_\ell\}$, draw $\ell$ independent functions $\psi_1,\ldots,\psi_\ell$, with $\psi_j$ following the distribution $P_{\bt_{\tau_j}}(\by_{\tau_j},df)$ conditioned on $f(\bt_{\tau_j^c})<\by_{\tau_j^c}$.
\item Independently of the above two steps, draw $\sum_{i\in I}\delta_{\phi_i}$ a Poisson point measure on $\bbC_0$ with intensity $1_{\{f(\bt)<\by\}}\mu(df)$. It can be obtained from a Poisson point measure with intensity $\mu(df)$ by removing those points not satisfying the constraint $f(\bt)<\by$.
\end{enumerate}
Then, the random field 
\[
\tilde \eta(t)=\max\{\psi_1(t),\ldots,\psi_\ell(t)\}\vee \max\{\phi_i(t),\ i\in I\},\quad t\in T,
\]
has the required conditional distribution. 

The issues and computational aspects of conditional sampling are addressed in the paper \cite{DEMR12}. The special case of Brown-Resnick max-stable processes is considered and  tractable expressions are derived and the above three-step procedure is implemented effectively.}\end{remark}

\section{Examples}
As an illustration, we apply in this section our general results to specific cases.

\subsection{The case of a single conditioning point}\label{sec:scp}
It is worth noting that  the case of a single conditioning point, i.e. $k=1$, gives rise to major simplifications. There exists indeed a unique partition of the set $K=\{t\}$ so that the notion of hitting scenario is irrelevant. Furthermore, there is a.s. a single $K$-extremal function $\varphi_1^+$ which is equal to the $t$-extremal function $\phi_{t}^+$. In this case, Theorems~\ref{theo2.2} and ~\ref{theo3} simplify into the following proposition.

\begin{proposition}\label{prop:scp} Let $t\in T$ and suppose that  conditions (i)-(iii) in Proposition~\ref{prop2} are met.
Then, conditionally on $\eta(t)=y$, $\phi_t^+$ and $\Phi_{\{ t \}}^-$ are independent; the conditional distribution of $\phi_t^+$ is equal to  $P_{t}(y,\cdot)$; the conditional distribution of $\Phi_{\{ t \}}^-$ is equal to the distribution of a Poisson point measure with intensity $1_{\{f(t)<y\}}\mu(df)$.
Furthermore, for $l\geq 1$, $\bs\in T^l$ and $\bz\in [0,+\infty)^l$,
\begin{eqnarray}
&&\bbP[\eta(\bs) < \bz  \mid \eta(t)=y]\label{eq:weintraub}\nonumber\\
&= &P_t(y,\{f(\bs)<\bz\})  \exp[-\mu(\{ f(\bs)\not<\bz, \  f(t)< y\})].
\end{eqnarray}
\end{proposition}

\subsection{Max-stable models}\label{sec:maxstable}
We put the emphasis here on max-stable random fields. For convenience
and without loss of generality, we focus on simple max-stable random
fields $\eta$, i.e., with standard unit Fr\'echet margins
\[
\bbP(\eta(t)\leq x)=\exp[-x^{-1}]1_{\{x>0\}},\quad x\in \bbR,\ t\in T.
\] 
A random field $\eta$ is said to be simple max-stable if for any $n\geq 1$, 
\[
 \eta \stackrel{\cL}= n^{-1}\bigvee_{i=1}^n \eta_i
\]
where $\{\eta_i,\ i\geq 1\}$ are i.i.d. copies of $\eta$. Any general max-stable random field can be related to
such a simple max-stable random field $\eta$ by simple transformation of the margins, see e.g. Corollary 3.6 in \cite{GHV90}.
Furthermore, Corollary 4.5.6 in \cite{dHF06} states that $\eta$ can be represented as
\begin{equation}\label{eq:repdH}
\{\eta(t)\}_{t\in T}\stackrel{\cL} =\Big\{\bigvee_{i\geq 1} \Gamma_i Y_i(t) \Big\}_{t\in T}
\end{equation}
where $(\Gamma_i)_{i\geq 1}$ is the nonincreasing enumeration of the points of a Poisson point process on $(0,\infty)$ with  intensity $x^{-2}dx$,
 $(Y_i)_{i\geq 1}$ is an i.i.d. sequence of  continuous random processes on $T$, independent of $(\Gamma_i)_{i\geq 1}$ and such that
\[
 \bbE[Y_1(t)]=1,\ t\in T, \quad \mbox{and} \quad  \bbE[\|Y_1\|]<\infty.
\]

Since a continuous simple max-stable random field is max-i.d., it has a Poisson point measure representation \eqref{repGHV}. The normalization to unit Fr\'echet margins entails that the vertex function $h$ is equal to $0$ and that the exponent measure $\mu$ satifies, for all $t\in T$,
\[
\mu_t(dy)=y^{-2}1_{\{y>0\}}dy\quad\mbox{and}\quad\bar\mu_t(y)=y^{-1},\quad
y>0.
\]
The correspondence between the two representations \eqref{repGHV} and \eqref{eq:repdH} is the following:  the point measure $\Phi=\sum_{i\geq 1} \delta_{\Gamma_i Y_i}$  is a Poisson point measure on $\bbC_0$ with intensity
\[
\mu(A)= \int_0^\infty\bbE[1_{\{rY_1\in A\}}]\, r^{-2}dr ,\quad A\subset \bbC_0 \mathrm{\ Borel\ set},
\]
The  distribution of the $Y_i$'s, denoted by $\sigma$, is called the spectral measure and is related to the exponent measure $\mu$ by the relation
\[
\mu(A)= \int_0^\infty\int_{\bbC_0} 1_{\{rf\in A\}}\, \sigma(df) r^{-2}dr ,\quad A\subset \bbC_0 \mathrm{\ Borel\ set}.
\]

Taking into account this particular form of the exponent measure, we can relate the kernel $P_t(y,\cdot)$ to the spectral measure $\sigma$. For $x\in\bbR$, we note $(x)^+=\max(x,0)$.
\begin{proposition}\label{prop4.2}
Let $\eta$ be a continuous simple max-stable random field with spectral measure $\sigma$ and $t\in T$. 
The $\{t\}$-extremal function  $\phi^+_t$  has conditional distribution
\[
\bbP[\phi_t^+\in \cdot \mid \eta(t)=y]=P_t(y,\cdot)=\int_{\bbC_0} 1_{\{ \frac{y}{f(t)}f\in\, \cdot\, \}}f(t)\sigma(df).
\]
Furthermore, for $l\geq 1$, $\bs\in T^l$ and $\bz\in [0,+\infty)^l$,
\begin{eqnarray}
&&\bbP[\eta(\bs) < \bz  \mid \eta(t)=y]\label{eq:weintraub2}\\
&=& \exp\Big[- \int_{\bbC_0}\Big( \vee_{i=1}^l \frac{f(s_{i})}{ z_{i} }-\frac{f(t)}{y}\Big)^+\sigma(df)\Big] \int_{\bbC_0}1_{\{ \vee_{i=1}^l   \frac{f(s_{i})}{z_{i}} <\frac{f(t)}{y} \}}f(t)\sigma(df) .\nonumber
\end{eqnarray}
\end{proposition}
Equation \eqref{eq:weintraub2} extends Lemma 3.4 in Weintraub \cite{W91} where only the  bivariate case $l=1$ is considered. Note the author considers  min-stability rather than max-stability; the correspondence is straightforward since, if $\eta$ is simple max-stable, then $\eta^{-1}$ is min-stable with exponential margins.

\subsection{Regular models}\label{sec:rm}
We have considered so far the case of a single conditioning point which allows for major simplifications.
In the general case, there are several conditioning points and the hitting scenario is non trivial. This introduces more complexity since the conditional distribution is expressed as a mixture over any possible hitting scenarios and involves an abstract Radon-Nykodym derivative. The framework of regular models can be helpful to get more tractable formulas.

The exponent measure $\mu$ is said to be {\it regular} (with respect
to the Lebesgue measure) if for any $l\geq 1$ and $\bs\in T^l$ with
pairwise distinct components, the measure $\mu_{\bs}(d\bz)$ is
absolutely continuous with respect to the Lebesgue measure $d\bz$ on
$[0,+\infty)^l$. We denote by $h_{\bs}$ the corresponding
Radon-Nykodym derivative, i.e., $\mu_{\bs}(d\bz)=h_{\bs}(\bz)d\bz$. 

Under this assumption, we can reformulate Theorems~\ref{theo2.1} and \ref{theo2.2}. For example, Equation \eqref{eq:theo2.1} implies that the distribution $\nu_\bt$ of $\eta(\bt)$ is  absolutely continuous with respect to the Lebesgue measure with density
\[
\frac{d\nu_{\bt}}{d\by}(\by)=\exp[-\bar\mu_\bt(\by)] \sum_{\tau\in\cP_K}\prod_{j=1}^{\ell(\tau)}\int_{\{\bz_j<\by_{\tau_j^c}\}} h_{(\bt_{\tau_j},\bt_{\tau_j^c})}(\by_{\tau_j},\bz_j)d\bz_j.  
\]
Equation \eqref{eq:theo2.2} giving the conditional distribution of the hitting scenario becomes
\[
\pi_\bt(\by,\tau)=\frac{\prod_{j=1}^{\ell(\tau)}  \int_{\{\bz_j<\by_{\tau_j^c}\}} h_{(\bt_{\tau_j},\bt_{\tau_j^c})}(\by_{\tau_j},\bz_j)d\bz_j}{\sum_{\tau'\in\cP_K}\prod_{j=1}^{\ell(\tau')}  \int_{\{\bz_j<\by_{{\tau_j'}^c}\}}h_{(\bt_{\tau_j'},\bt_{{\tau_j'}^c})}(\by_{\tau_j'},\bz_j)d\bz_j}.
\]
The conditional distribution of the extremal functions $Q_\bt(\by,\tau,\cdot)$ in Equation \eqref{eq:theo2.3} is based on the kernel $P_\bt(\by,df)$. Using the existence of a Radon-Nykodym derivative for the finite dimensional margins of $\mu$, we obtain 
\[
P_\bt(\by, f(\bs)\in d\bz)=\frac{h_{(\bt,\bs)}(\by,\bz)}{h_{\bt}(\by)}d\bz.
\]
This approach is exploited in \cite{DEMR12} for  Brown-Resnick max-stable processes. Indeed, the model turns out to be regular.

\section{Proofs}
\subsection{Proof of Theorem~ \ref{theo1} and Proposition \ref{prop1}}
For the proof of Theorem \ref{theo1}, we need the following lemma giving a useful characterization of the $K$-extremal random point measure. If $M_1,M_2\in M_p(\bbC_0)$ are such that 
$M_2-M_1\in M_p(\bbC_0)$, we call $M_1$  a sub-point measure of $M_2$.
\begin{lemma}\label{lem1}
The $K$-extremal  point measure $\Phi_K^+$ is the   unique sub-point measure $\tilde \Phi$ of $\Phi$ such that
\[
\tilde \Phi \in C^+_K \quad \mathrm{and}\quad  \Phi-\tilde \Phi\in C^-_K(\max(\tilde \Phi)).
\]
\end{lemma}
\paragraph {Proof of Lemma \ref{lem1}:} 

 First the condition $\Phi-\tilde \Phi\in
  C^-_K(\max(\tilde \Phi))$ implies
\[
\max(\Phi-\tilde \Phi)<_K \max(\tilde \Phi)\quad \mbox{and}\quad \max(\tilde \Phi)=_K\max(\Phi).
\]
Let $f\in [\Phi-\tilde \Phi]$. The condition $\Phi-\tilde \Phi\in C^-_K(\max(\tilde \Phi))$ implies  $f<_K \max(\tilde \Phi)$. Since $\tilde \Phi$ is a sub-point measure of $\Phi$, $\max(\tilde \Phi)\leq\max(\Phi)$ so that   $f<_K \max(\Phi)$ and  $f$ is $K$-sub-extremal in $\Phi$.\\
Conversely for $f\in [\tilde \Phi]$, the condition $\tilde \Phi\in
C^+_K$ implies the existence of $t_0\in K$
such that $f(t_0)=\max(\tilde \Phi)(t_0)$.  Hence $f(t_0)=\max(\Phi)(t_0)$
and $f$ is $K$-extremal in $\Phi$.\\
\qed

\paragraph {Proof of Theorem \ref{theo1} :}  
  First note that $\Phi_K^+(\bbC_0)=0$ if and only if $\Phi=0$. This occurs with probability $\exp[-\mu(\bbC_0)]$ and in this case $\Phi_K^+=\Phi_K^-=0$. The first claim follows.\\
  Next, let $k\geq 1$. According to Lemma \ref{lem1},
  $\Phi_K^+(\bbC_0)=k$ if and only if there exists a $k$-uplet
  $(\phi_1,\ldots,\phi_k)\in[\Phi]^{k}$ such that
\[
\sum_{i=1}^k \delta_{\phi_i} \in C_K^+ \quad \mathrm{and} \quad \Phi-\sum_{i=1}^k\delta_{\phi_i}\in C_K^-\left(\vee_{i=1}^k \phi_i \right) .
\]
When this holds, the $k$-uplet $(\phi_1,\ldots,\phi_k)$ is unique up
to a permutation of the coordinates and we have
\[
 \sum_{i=1}^k \delta_{\phi_i}= \Phi_K^+ \quad \mathrm{and} \quad \Phi-\sum_{i=1}^k\delta_{\phi_i}= \Phi_K^-.
\]
Hence the sum
\begin{eqnarray*}
&&\int_{\bbC_0^k} 1_{ \left\lbrace  \sum_{i=1}^k \delta_{\phi_i}\in A\cap  C_K^+,\ \Phi-\sum_{i=1}^k\delta_{\phi_i}\in C_K^- \left(  \vee_{i=1}^k \phi_i \right)  \right\rbrace }\,\Phi(d\phi_1)\cdots\big( \Phi-\sum_{j=1}^{k-1} \delta_{\phi_j} \big) (d\phi_k)
\end{eqnarray*}
is equal to 
$k! 1_{\{(\Phi_K^+,\Phi_K^-)\in A\times B \}}$ if $\Phi_K^+(\bbC_0)=k$ and $0$ otherwise.
Using this and Slyvniak's formula (see Appendix~\ref{sec:slyvniak}), we get
\begin{eqnarray*} 
&&k!\bbP\big[(\Phi_K^+,\Phi_K^-)\in A\times B  ,\ \Phi_K^+(\bbC_0)=k\big]\\
&=&  \bbE \Big[\int_{\bbC_0^k} 1_{ \left\lbrace  \sum_{i=1}^k \delta_{\phi_i}\in A\cap  C_K^+,\ \Phi-\sum_{i=1}^k\delta_{\phi_i}\in B\cap C_K^-\left(\vee_{i=1}^k \phi_i \right)  \right\rbrace }\\
&&\qquad\Phi(d\phi_1)\cdots \big( \Phi-\sum_{j=1}^{k-1} \delta_{\phi_j} \big) (d\phi_k)\Big]\\
&=&\int_{\bbC_0^k}1_{\{\sum_{i=1}^k \delta_{f_i}\in A\cap C^+_K\}} \bbP \left[\Phi \in B \cap C^-_K \left(\vee^k_{i=1} f_i \right) \right]\,\mu^{\otimes k}(df_1,\ldots,df_k).
\end{eqnarray*}
This proves Theorem \ref{theo1}.
\qed

 \paragraph {Proof of Proposition \ref{prop1}:} 
  In the case $\mu(\bbC_0)<+\infty$, $\Phi$ and {\it a fortiori}
  $\Phi_K^+$ are a.s. finite.  Suppose now
  $\mu(\bbC_0)=+\infty$, so that $\Phi$ is a.s. infinite. If
  $\inf_{t\in K} \eta(t)=0$, then there is $t_0\in K$ such that
  $\eta(t_0)=0$ (recall $\eta$ is continuous and $K$ compact). This
  implies that $\phi(t_0)=0$ for all $\phi\in[\Phi]$ and hence
  $\Phi_K^+=\Phi$ is infinite. If $\inf_{t\in K}
  \eta(t)=\varepsilon >0$, then the support of $\Phi_K^+$ is
  included in the set $\{f\in\bbC_0;\ f(K)\geq \varepsilon\}$. From
  the definition of $M_p(\bbC_0)$, this set contains only a finite
  number of atoms of $\Phi$ so that $\Phi_K^+$ must be finite.
\qed

\subsection{Proof of Propositions \ref{prop2} and \ref{prop3}}

\paragraph {Proof of Proposition~\ref{prop2}:} 
According to equation \eqref{eq:C2bis}, for all $x>0$, 
\[
 \bbP[\eta(t)= x] = \exp[-\bar\mu_t(x^{+})] -\exp[-\bar\mu_t(x)],
\]
and $ \bbP[\eta(t)= 0]=\exp[-\bar\mu_t(0^{+})]$. The equivalence between (ii) and (iii) follows.\\
The equivalence between (i) and (ii) is a consequence of
Theorem \ref{theo1} with $K=\{t\}$, $k=1$ and $A=B=M_p(\bbC_0)$: we get
\begin{eqnarray*}
\bbP[\Phi_{\{t\}}^+(\bbC_0)=1]&=& \int_{\bbC_0} 1_{\{\delta_f\in C_{\{t\}}^+ \}}\bbP[ \Phi\in C_{\{t\}}^-(f)]\mu(df)\\
&=&\int_{[0,+\infty)}\exp[-\bar\mu_t(y)]\,\mu_t(dy).
\end{eqnarray*} 
It remains to prove that this probability is equal to $1$ if and only
if (ii) is satisfied. To this aim, we compute
\begin{eqnarray}
\bbP[\Phi_{\{t\}}^+(\bbC_0)=1]\label{eq:intcond}
&=&\int_{(0,+\infty)^{2}} e^{-x}1_{\{ x \geq \bar\mu_t(y ) \}}\,dx\mu_t(dy) \\
&=&\int_{ (0,+\infty) }e^{-x}\mu_t(A_{x})\,dx,\nonumber 
\end{eqnarray}
where $A_{x}=\{y> 0:   \bar\mu_t(y ) \leq x\}$. Since $\bar
\mu_t$ is c\`ag-l\`ad, non-increasing and tends to $\infty$ at $0$, 
$A_x = (\inf A_x, \infty ) \neq \varnothing$ for all $x
> 0$. Furthermore using equation~\eqref{eq:intcond} and the fact that
$\mu_t(A_x) \leq x$,  
we get that $\bbP[\Phi_{\{t\}}^+(\bbC_0)=1]=1$ if and
only if $\mu_t(A_x)=x$ for all $x>0$.  We see easily that this is 
equivalent to condition (ii) and this completes the equivalence between (i)
and (ii).

We now prove Equation \eqref{eq:prop2.2}. Assuming that conditions (i)-(iii) are met, it holds
\[
\bbP(\phi_t^+\in A)=\bbP[\Phi_{\{t\}}^+\in \tilde A,\  \Phi_{\{t\}}^+(\bbC_0)=1]
\]
with $\tilde A=\{\delta_f,\ f\in A\}$. Theorem \ref{theo1} with $K=\{t\}$, $k=1$ and $B=M_p(\bbC_0)$ entails
\begin{eqnarray*}
\bbP(\phi_t^+\in A)&=&\int_{\bbC_0}1_{\{ \delta_{f}\in \tilde A\cap C^+_{\{t\}}\}} \bbP[\Phi\cap C^-_{\{t\}}(f)] \,\mu(df)\\
&=&\int_{\bbC_0}1_{\{f\in A \}} \exp[-\bar\mu_{t}(f(t))]\,\mu(df).
\end{eqnarray*}
This proves Equation \eqref{eq:prop2.2}.
\qed

 \paragraph {Proof of Proposition \ref{prop3}:}
First note that the inequalities \eqref{eq:ineqhs} characterize the hitting scenario. Let $\tau=(\tau_1,\ldots,\tau_\ell)\in\cP_K$  and define the sets
\[
\tilde C_\tau=\Big\{  (f_{1},\ldots,f_{\ell})\in \bbC_0^\ell;\ 
\forall j\in[\![1,\ell]\!],\ f_j>_{\tau_j} \bigvee_{j'\neq j} f_{j'} \Big\}.
\]
and
\[
C_\tau=\Big\{ \sum_{j=1}^\ell \delta_{f_j}\in M_p(\bbC_0);\  (f_{1},\ldots,f_{\ell})\in \tilde C_\tau \Big\}.
\]
Note that $C_\tau\subset C_K^+$ and that $\Theta=\tau$ if and only if $\Phi_{K}^+\in C_\tau$. \\
Furthermore, $\Theta=\tau$ and $(\varphi^+_{1},\ldots,\varphi^+_{\ell})\in A$ if and only if $\Phi_K^+\in A_\tau$ with
\[
A_\tau=\Big\{ \sum_{j=1}^\ell \delta_{f_j}\in M_p(\bbC_0);\  (f_{1},\ldots,f_{\ell})\in   C_\tau\cap A \Big\}.
\]
Hence the following events are equal
\[
\{ \Theta=\tau,\ (\varphi^+_1,\ldots,\varphi^+_\ell)\in A,\ \Phi_K^-\in B\}
=\{ \Phi_K^+ \in  A_\tau,\ \Phi_K^-\in B,\ \Phi_K^+(\bbC_0)=\ell\}
\]
and Theorem \ref{theo1} implies
\begin{eqnarray}
&&\bbP[\Theta=\tau,\ (\varphi^+_1,\ldots,\varphi^+_\ell)\in A,\ \Phi_K^-\in B]\nonumber\\
&=& \frac{1}{\ell !}\int_{\bbC_0^\ell}1_{\{\sum_{j=1}^\ell \delta_{f_j}\in A_\tau\}}
\bbP \left[\Phi \in B \cap C^-_K \left(\vee^\ell_{j=1} f_j \right) \right]
\,\mu^{\otimes \ell}(df_1,\ldots,df_\ell)\label{eq:prop3.2}.
\end{eqnarray}
Finally, $\sum_{j=1}^\ell \delta_{f_j}\in  A_\tau$ if and only
if there exists a permutation $\sigma$ of $[\![1,\ell]\!]$ such that
$(f_{\sigma(1)},\ldots,f_{\sigma(\ell)}) \in A\cap \tilde C_\tau$. Such a
permutation is unique and this proves the equivalence of Equations \eqref{eq:prop3.1} and \eqref{eq:prop3.2}.
\qed

\subsection{Proofs of Theorems \ref{theo2.1}, \ref{theo2.2} and \ref{theo3}}
\paragraph {Proof of Theorems \ref{theo2.1} and \ref{theo2.2}:} 
Note that $\eta(\bt)$ can be expressed in terms of the hitting
scenario and the extremal function as follows. For $\tau\in\cP_K$,
define the mapping $\Gamma_\tau:\bbC_0^\ell \to [0,+\infty)^k$ by
\[
\Gamma_\tau(f_1,\ldots,f_\ell)=(y_1,\ldots,y_k)\quad \mathrm{ with}\ y_i=f_j(t_i)\ \mathrm{ if}\ t_i\in \tau_j.
\]
Definition \eqref{def:hs} entails that for all $t\in\theta_j$, $\eta(t)=\varphi_j^+(t)$. This can be rewritten as  $\eta(\bt)=\Gamma_{\Theta}(\phi_{1}^+,\ldots,\phi_{\ell}^+)$.
Using this, the probability
\[
P(\tau,A,B,C)=\bbP[\Theta=\tau,\ (\varphi^+_{_1},\ldots,\varphi^+_{_\ell})\in A,\ \Phi_K^-\in B,\ \eta(\bt)\in C]
\]
can be computed thanks to Proposition \ref{prop3}:
 \begin{eqnarray*}
&&P(\tau,A,B,C)\\
&=& \bbP[\Theta=\tau,\  (\varphi^+_{_1},\ldots,\varphi^+_{_\ell})\in A\cap \Gamma_\tau^{-1}(C),\ \Phi_K^-\in B,\ \eta(\bt)\in C]\\
&=& \int_{A\cap \Gamma_\tau^{-1}(C)}1_{ \left\{\forall j\in[\![1,\ell]\!],\ f_j>_{\tau_j} \bigvee_{j'\neq j} f_{j'} \right\}}
\bbP \left[\Phi \in B \cap C^-_K \left(\vee^\ell_{j=1} f_j \right) \right]\,\mu^{\otimes \ell}(df_1,\ldots,df_\ell). 
\end{eqnarray*}
Now for each $j\in [\![1,\ell]\!]$, we condition  the measure $\mu(df_j)$ with respect to $f_j(\bt_{\tau_j})$: Equation \eqref{eq:theo2.1pre} entails
\begin{eqnarray}
&&P(\tau,A,B,C) \label{eq:theo2.1bis}\\
&=&\int_{C}\int_{A} 1_{\left\{\forall j\in[\![1,\ell]\!],\ f_j>_{\tau_j} \bigvee_{j'\neq j} f_{j'} \right\}} \bbP\left[\Phi \in B \cap C^-_K \left(\vee^\ell_{j=1} f_j \right) \right]\nonumber \\
&&\qquad\qquad \otimes _{j=1}^\ell  \big\{P_{\bt_{\tau_j}}(\by_{\tau_j},df_j)\, \mu_{\bt_{\tau_j}}(d\by_{\tau_j})\big\}\nonumber\\
&=& \int_{C}\int_{A} 1_{\{\forall j\in[\![1,\ell]\!],\ f_j(\bt_{\tau_j^c}) <\, \by_{\tau_j^c}\}}\bbP[\Phi \in B \cap C^-_\bt(\by) ]\nonumber\\
&&\qquad\qquad \otimes _{j=1}^\ell  \big\{P_{\bt_{\tau_j}}(\by_{\tau_j},df_j)\, \mu_{\bt_{\tau_j}}(d\by_{\tau_j})\big\}.\nonumber
\end{eqnarray}
In the last equality, we use the fact that $f_j(\bt_{\tau_j})=\by_{\tau_j}$ a.s. under $ P_{\bt_{\tau_j}}(\by_{\tau_j},df_j)$, whence
\begin{eqnarray*}
\left\{\forall j\in[\![1,\ell]\!],\ f_j>_{\tau_j} \vee_{j'\neq j} f_{j'}\right\} 
&=&\left\{\forall j\in[\![1,\ell]\!],\  f_j<_{\tau_j^c} \vee_{j'\neq j} f_{j'}\right\}\\
&=&\left\{\forall j\in[\![1,\ell]\!],\ f_j(\bt_{\tau_j^c}) < \by_{\tau_j^c}\right\} 
\end{eqnarray*}
and $C^-_K (\vee^\ell_{j=1} f_j )=C^-_\bt(\by)$.
 
We now prove Theorem \ref{theo2.1}. Setting $A=\bbC_0^\ell$ and $B=M_p(\bbC_0)$ in Equation \eqref{eq:theo2.1bis}, we obtain
\begin{eqnarray*}
&&\bbP[\Theta=\tau,\ \eta(\bt)\in C]\\
&=&\int_{C} \int_{\bbC_0^\ell} 1_{\{\forall j\in[\![1,\ell]\!],\ f_j(\bt_{\tau_j^c}) <\, \by_{\tau_j^c}\}}\bbP[\Phi \in  C^-_\bt(\by) ]\, \otimes _{j=1}^\ell  \big\{P_{\bt_{\tau_j}}(\by_{\tau_j},df_j)\, \mu_{\bt_{\tau_j}}(d\by_{\tau_j})\big\}.
\end{eqnarray*}
Using the fact that $\bbP[\Phi \in C^-_\bt(\by)
]=\exp[-\bar\mu_{\bt}(\by)]$ and performing integration with respect
to $\otimes_{j=1}^\ell P_{\bt_{\tau_j}}(\by_{\tau_j},df_j)$, we obtain Equation \eqref{eq:theo2.1} and this proves Theorem \ref{theo2.1}. 

We now consider Theorem \ref{theo2.2}. Combining Equations \eqref{eq:theo2.1}-\eqref{eq:theo2.4} together with Equation \eqref{eq:theo2.1bis}, we get
\begin{eqnarray}
&& \bbP[\Theta=\tau,\  (\varphi^+_{_1},\ldots,\varphi^+_{_\ell})\in A,\ \Phi_K^-\in B,\ \eta(\bt)\in C]\nonumber\\
&=& \int_{C}\int_{A}  \frac{\bbP[\Phi \in B \cap C^-_\bt(\by) ]}{\bbP[\Phi \in  C^-_\bt(\by)]}\, Q_\bt(\by,\tau,df_1\cdots df_\ell) \nu_\bt^\tau(d\by)\nonumber\\
&=& \int_{C}\int_{A} R_\bt(\by,\tau,(f_j),B) Q_\bt(\by,\tau,df_1\cdots df_\ell) \pi_\bt(\by,\tau) \nu_\bt(d\by).\label{eq:theo2.1ter}
\end{eqnarray}
In particular, with $A=\bbC_0^\ell$ and $B=M_p(\bbC_0)$, we obtain the relation
\[
\bbP[\Theta=\tau,\ \eta(\bt)\in C]= \int_{C}\pi_\bt(\by,\tau)\nu_\bt(d\by) 
\]
characterizing the regular conditional distribution $\bbP[\Theta=\tau\mid \eta(\bt)=\by]$ (see Appendix~\ref{sec:rcd}). This proves that Equation \eqref{eq:theo2.2} provides the regular conditional distribution $\bbP[\Theta=\tau\mid \eta(\bt)=\by]$.
Similarly, Equation \eqref{eq:theo2.1ter} entails that the regular conditional distributions 
$\bbP[(\varphi_j^+)\in \cdot \mid \eta(\bt)=\by, \Theta=\tau]$ and $\bbP[\Phi_K^-\in \cdot \mid \eta(\bt)=\by, \Theta=\tau, (\varphi_j^+)=(f_j)]$ are given respectively by $Q_\bt(\by,\tau,\cdot)$ in Equation \eqref{eq:theo2.3} and  $R_\bt(\by,\tau,(f_j),\cdot)$ in Equation \eqref{eq:theo2.4}. 

We briefly comment on these formulas. The fact that the distribution $Q_\bt(\by,\tau,\cdot)$ in Equation \eqref{eq:theo2.3} factorizes into a tensorial product means that the extremal functions $\varphi_1^+,\ldots,\varphi_\ell^+$ are independent conditionally on $\eta(\bt)=\by$ and $\Theta=\tau$. The fact that the distribution $R_\bt(\by,\tau,(f_j),\cdot)$  in Equation \eqref{eq:theo2.4} does not depend on $\tau$ and $(f_j)$ means that conditionally on $\eta(\bt)=\by$, $\Phi_K^-$ is independent of $\Theta$ and $(\varphi_1^+,\ldots,\varphi_{\ell(\Theta)}^+)$. The  distribution $R_\bt(\by,\cdot)$ can be seen as the distribution of the Poisson point measure $\Phi$ conditioned to lie in $C_\bt^-(\by)$, i.e., to have no atom in $\{f\in\bbC_0;\ f(\bt)\not < \by\}$. It is equal to the distribution of a Poisson point measure with intensity $1_{\{f(\bt)<\by\}}\mu(df)$. 
\qed

\paragraph {Proof of Theorem~\ref{theo3}:} 
Remark that
\begin{eqnarray*}
\{\eta(\bs)< \bz \}&=&\{(\Phi_K^+,\Phi_K^-)\in C_\bs^-(\bz)\times C_\bs^-(\bz) \}\\
&=& \bigcup_{\tau\in\cP_K} \{\Theta=\tau, \varphi_1^+(\bs)<\bz,\ldots,\varphi_\ell^+(\bs)<\bz, \Phi_K^-\in C_\bs^-(\bz) \}
\end{eqnarray*}
where $C_\bs^-(\bz)$ is defined in Theorem \ref{theo2.2}. Using this, Theorem \ref{theo2.2} entails
\begin{eqnarray*}
&&\bbP[\eta(\bs)< \bz \mid \eta(\bt)=\by]\\
&=&\sum_{\tau\in\cP_K} \bbP[\Theta=\tau, \varphi_1^+(\bs)<\bz,\ldots,\varphi_\ell^+(\bs)<\bz, \Phi_K^-\in C_\bs^-(\bz) \mid \eta(\bt)=\by]\\
&=& \sum_{\tau\in\cP_K}\pi_\bt(\by,\tau)Q_\bt(\by,\tau,\{f(\bs)<\bz\}^\ell)R_\bt(\by, C_\bs^-(\bz)).
\end{eqnarray*}
The result follows since
\[
Q_\bt(\by,\tau,\{f(\bs)<\bz\}^\ell) = \prod_{j=1}^\ell \frac{P_{\bt_{\tau_j}}( \by_{\tau_j},\{ f(\bt_{\tau_j^c})< \by_{\tau_j^c}, \ f(\bs)<\bz\})}{P_{\bt_{\tau_j}}(\by_{\tau_j},\{ f(\bt_{\tau_j^c})< \by_{\tau_j^c}\} )}
\]
and
\begin{eqnarray*}
R_\bt(\by, C_\bs^-(bz))&=&\frac{\bbP[\Phi \in C^-_\bs(\bz) \cap C^-_\bt(\by) ]}{\bbP[\Phi \in  C^-_\bt(\by) ]}\\
&=&\frac{\exp[-\mu(\{ f(\bs)\not<\bz \ \mathrm{or}\ f(\bt)\not< \by\})]}{\exp[-\mu(\{  f(\bt)\not< \by\})]}\\
&=&\exp[-\mu(\{ f(\bs)\not<\bz, \ f(\bt)< \by\})].
\end{eqnarray*}
\qed

\subsection{Proof of Propositions \ref{prop:scp} and \ref{prop4.2}}
\paragraph {Proof of Proposition \ref{prop:scp}:} 
This is a straightforward application of Theorem \ref{theo2.2} and \ref{theo3}. Take into account that when $K=\{t\}$, $\cP_K$ is reduced to a unique partition of size $\ell=1$ so that $\Theta=\{t\}$  and  $\varphi_1^+=\phi_t^+$.
\qed

\paragraph {Proof of Proposition \ref{prop4.2}:} 
According to Proposition \ref{prop:scp}, $\bbP[\phi_t^+\in \cdot \mid \eta(t)=y]$ is equal to$P_t(y,\cdot)$. 
For any measurable $A\subset\bbC_0$ and $B\subset (0,+\infty)$, we compute 
\begin{eqnarray*}
&&\int_B \int_{\bbC_0} 1_{\{ \frac{y}{f(t)}f\in A \}}f(t)\sigma(df) \mu_t(dy)\\
&=& \int_0^\infty \int_{\bbC_0}\int_{\bbC_0}1_{\{ \frac{rg(t)}{f(t)}f\in A \}}1_{\{rg (t)\in B\}}f(t)\sigma(df) r^{-2}dr\sigma(dg)\\
&=& \int_{\bbC_0}  \int_0^\infty 1_{\{ r f\in A  \}}  1_{\{rf (t)\in B\}}\, r^{-2}dr \sigma(df)\\
&=&  \int_{\bbC_0}1_{\{f\in A\}}1_{\{f(t)\in B\}}\,\mu(df)
\end{eqnarray*}
The second equality follows from the change of variable $\tilde r = rg(t)/f(t)$ together with  the relation $\int_{\bbC_0} g(t)\sigma(dg)=1$.
This proves that
\[
P_t(y,A)=\int_{\bbC_0} 1_{\{ \frac{y}{f(t)}f\in A \}}f(t)\sigma(df).
\]
According to Equation~\eqref{eq:weintraub}
\[
\bbP[\eta(\bs) < \bz  \mid \eta(t)=y]= P_t(y,\{f(\bs)<\bz\})  \exp[-\mu(\{ f(\bs)\not<\bz, \  f(t)< y\})]. 
\]
We have 
\begin{eqnarray*}
P_t(y, \{f(\bs) < \bz \} )&=& \int_{\bbC_0} 1_{\{ \frac{y}{f(t)}f(\bs)<\bz \}}f(t)\sigma(df)\\
&=&\int_{\bbC_0}1_{\{ \vee_{i=1}^l  \frac{f(s_{i})}{ z_{i} } <\frac{f(t)}{y} \}}f(t)\sigma(df)
\end{eqnarray*}
and 
\begin{eqnarray*}
\mu(\{ f(\bs) \not< \bz,\, f(t) < y\} )&=&\int_0^\infty \int_{\bbC_0} 1_{\{ rf(\bs) \not< \bz,\, rf(t) < y\}}r^{-2}dr \sigma(df)\\ 
&=& \int_0^\infty \int_{\bbC_0} 1_{\{    \min\limits_{1\leq i\leq l}  \frac{z_{i}}{f(s_{i})} \leq r < \frac{y}{ f(t)} \}}r^{-2}dr \sigma(df) \\ 
&=& \int_{\bbC_0}   \Big( \vee_{i=1}^l  \frac{f(s_{i})}{ z_{i} } - \frac{f(t)}{y}  \Big)^+    \sigma(df).
\end{eqnarray*}
This proves Equation~\eqref{eq:weintraub2}.
\qed

\appendix
\section{Auxiliary results}
\subsection{Slyvniak's formula}\label{sec:slyvniak}

Palm Theory deals with conditional distribution for point
processes. We recall here one of the most famous formula of Palm
theory, known as Slyvniak's Theorem. This will be the main tool in our
computations. For a general reference on Poisson point processes, Palm
theory and their applications, the reader is invited to refer to the
monograph \cite{SKM87} by Stoyan, Kendall and Mecke.

Let $M_{p}(\bbC_0)$ be the set of locally-finite point measures $N$ on $\bbC_0$ endowed with the $\sigma$-algebra generated by the family of mappings $$N\mapsto N(A),\qquad  A\subset \bbC_0 \mathrm{\ Borel\ set}.$$ 

\begin{theorem}[Slyvniak's Formula]\ \\
Let $\Phi$ be a Poisson point process on $\bbC_0$ with intensity measure $\mu$. For any measurable function $F:\bbC_0^k\times M_{p}(\bbC_0)\to [0,+\infty)$, 
\begin{eqnarray*}
& &\bbE\Big[\int_{\bbC_0^k} 
  F\Big(\phi_1,\ldots,\phi_k,\Phi-\sum_{i=1}^k \delta_{\phi_i}\Big)\,\Phi(d\phi_1)\,(\Phi-\delta_{\phi_1})(d\phi_2)\cdots \Big(\Phi-\sum_{j=1}^{k-1} \delta_{\phi_j}\Big)(d\phi_k) \Big]\\
 &=&\int_{\bbC_0^k}\bbE[F(f_1,\ldots,f_k,\Phi)]\,\mu^{\otimes k}(df_1,\ldots,df_k).
\end{eqnarray*}
\end{theorem}

\subsection{Regular conditional distribution}\label{sec:rcd}
We recall here briefly the notion of regular conditional probability (see e.g. Proposition A1.5.III in Daley and Vere-Jones \cite{DVJ03}).
Let $(\cY,\cG)$ be a complete separable metric space with its associated $\sigma$-algebra of Borel sets, $(\cX,\cF)$ an arbitrary measurable space, and $\pi$ a probability measure on the product space $(\cX\times\cY,\cF\otimes\cG)$. Let $\pi_\cX$ denote the $\cX$-marginal of $\pi$, i.e.  $\pi_\cX(A)=\pi(A\times\cY)$ for any $A\in\cF$. Then  there exists a family of kernels $K(x,B)$ such that 
\begin{itemize}
\item[-] $K(x,\cdot)$ is a probability measure on $(\cY,\cG)$ for any fixed $x\in\cX$;
\item[-] $K(\cdot,B)$ is an $\cF$-measurable function on $\cX$ for each fixed $B\in\cG$;
\item[-] $\pi(A\times B)=\int_A K(x,B)\pi_\cX(dx)$ for any $A\in\cF$ and $B\in\cG$.
\end{itemize}
These three properties define the notion of regular conditional probability. 
When $\pi$ is the joint distribution of the random variable $(X,Y)$, we may write $K(x, \cdot )  = \mathbb P (Y \in \cdot | X=x)$.

The existence of the regular conditional probability relies on  the assumption that $\cY$ is a complete and separable metric space. Furthermore, 
for any $\cF\otimes\cG$-measurable non-negative function $f$ on $X\times Y$, it follows that
\[
 \int_{\cX\times \cY}f(x,y)\pi(dx,dy)=\int_{\cX}\int_{\cY}f(x,y)K(x,dy)\pi_\cX(dx).
\]

The following Lemma states the existence of the kernel $$\{P_{\bs}(\bx,df);\ \bx\in[0,+\infty)^l\setminus\{0\}\}$$ satisfying Equation~\eqref{eq:theo2.1pre}. This is not straightforward since the measure $\mu$ is not a probability measure and may be infinite. 
\begin{lemma}\label{lem:rcd}
The regular version of the conditional measure $\mu(df)$ with respect to $f(\bs)\in [0,+\infty)^l\setminus\{0\}$ exists. It is denoted by $\{P_{\bs}(\bx,df);\ \bx\in[0,+\infty)^l\setminus\{0\}\}$  and satisfies Equation~\eqref{eq:theo2.1pre}.
\end{lemma}
\paragraph {Proof :}
Let  $|\cdot|$ denote a norm on $[0,+\infty)^l$. Define $A=\{f\in\bbC_0;\ f(\bs)\neq 0\}$ and, for $i\geq 0$, $A_i=\{f\in\bbC_0;\ (i+1)^{-1}\leq |f(\bs)|< i^{-1}\}$ with the convention $0^{-1}=+\infty$. Clearly, $A$ is equal to the disjoint union of the $A_i$'s. We note $\mu^i(\cdot)=\mu(\cdot\cap A_i)$ the measure on the complete and separable space $\bbC_0\cup\{0\}$. Equation \eqref{eq:condmu} ensures that $\mu_i$ is a finite measure (and hence a probability measure up to a normalization constant) and  there exists a regular conditional probability kernel $P^i_\bs(\bx,df)$ with respect to $f(\bs)=\bx$. We obtain, for all $F: [0,+\infty)^{l}\times\bbC_0$,
\[
 \int_{A_i}F(f(\bs),f)\mu(df)=\int_{\tilde A_i}\int_{\bbC_0}F(\bx,f)P^i_\bs(\bx,df)\mu_\bs(d\bx),
\]
where $\tilde A_i=\{x\in [0,+\infty)^k;\ (i+1)^{-1}\leq |\bx|< i^{-1}\}$. Let us define 
$P_\bs(\bx,df)$ a probability  measure on $\bbC_0$ by 
$$
P_\bs(\bx,df)= \sum_{i\geq 1} 1_{ \{ \bx\in\tilde A_i \} } P^i_\bs(\bx,df).
$$
If $F$ vanishes on $\{0\}\times\bbC_0$, we obtain
\begin{eqnarray*}
 \int_{\bbC_0}F(f(\bs),f)\mu(df)&=&\sum_{i\geq 0} \int_{A_i}F(f(\bs),f)\mu(df)\\
&=& \sum_{i\geq 0} \int_{\tilde A_i}\int_{\bbC_0}F(\bx,f)P^i_\bs(\bx,df)\mu_\bs(d\bx)\\
&=& \int_{[0,+\infty)^k\setminus\{0\} }\int_{\bbC_0}F(\bx,f)P_\bs(\bx,df)\mu_\bs(d\bx).
\end{eqnarray*}
This proves Equation~\eqref{eq:theo2.1pre}.
\qed

\subsection{Measurability properties}\label{sec:meas}
\begin{lemma}\label{lem:meas1}
 $\Phi_K^+$ and $\Phi_K^-$ are measurable from $(\Omega,\cF,\bbP)$ to $(M_p(\bbC_0),\cM_p)$.
\end{lemma}
\paragraph {Proof :}
From  Definition \ref{def:ep}, it is enough to prove that that the events $\{\phi_i<_K\eta\}\in\cF$ and $\{\phi_i\not<_K\eta\}$ are $\cF$-measurable. Let $K_0$ be a dense countable subset of $K$ and note that $\phi<_K\eta$ if and only if there is some rational $\varepsilon>0$ so that $\phi(t)<\eta(t)-\varepsilon$ for all $t\in K_0$.  Hence, for all $ n\in \bbN\cup\{+\infty\}$,
\begin{eqnarray}
 \{\phi_i<_K \eta\ ;\ N=n\}&=& \bigcup_{\varepsilon >0}\bigcap_{t\in K_0}  \{N=n\ ;\ \phi_i(t) <\eta(t)-\varepsilon\}\nonumber\\
&=& \bigcup_{\varepsilon >0}\bigcap_{t\in K_0} \bigcup_{j\leq n} \{N=n\ ;\ \phi_i(t)<\phi_j(t)-\varepsilon\}\label{eq:ssext}
\end{eqnarray}
and $\{\phi_i<_K \eta\}=\bigcup_{n=0}^\infty \{\phi_i<_K \eta\ ;\ N=n\}\in \cF$. Note the union over $\varepsilon$ is countable since $\varepsilon$ is taken rational.
\qed
\begin{lemma}\label{lem:meas2}
The set $C^+_K$ and $C^-_K(g)$ are measurable in $(M_p(\bbC_0),\cM_p)$. 
\end{lemma}
\paragraph {Proof :}
Let $g$ be a continuous function defined at least on $K$ and consider the Borel set  
\[
A=\{f\in\bbC_0;\ f\not<_K g\}\subset\bbC_0.
\]
The set $C^-_K(g)$ defined by Equation \eqref{eq:Cflat} is equal to  
\[
C^-_K(g)=\{M\in M_p(\bbC_0);\ M(A)=0\}
\]
and is $\cM_p$-measurable.\\
In order to prove the measurability of $C_K^+$ defined by Equation \eqref{eq:Csharp}, we introduce a measurable enumeration of the atoms of a point measure $M$ (see Lemma 9.1.XIII in  Daley and Vere-Jones \cite{DVJ08}). 
One can construct measurable applications 
$$\kappa:M_p(\bbC_0)\to \bbN\cup\{\infty\}\quad \mbox{ and}\quad \psi_i:M_p(\bbC_0)\to\bbC_0,\ i\geq 1,$$
such that 
\[
M=\sum_{i=1}^{\kappa(M)}\delta_{\psi_i(M)},\quad M\in M_p(\bbC_0).
\] 
A point measure $M$ does not lie in $C_K^+$ if and only if it has a $K$-subsextremal atom. Hence, 
\[
M_p(\bbC_0)\setminus C^+_K= \bigcup_{k=0}^{+\infty}\bigcup_{i=1}^k \{\kappa(M)=k;\ \psi_i(M)<_K\max(M) \}.
\]
Similar computations as in Equation \eqref{eq:ssext} entail
\[
\{\kappa(M)=k;\ \psi_i(M)<_K\max(M) \}= \bigcup_{\varepsilon >0}\bigcap_{t\in K_0} \bigcup_{j\leq k} \{\kappa=k\ ;\ \psi_i(t) <\psi_j(t)-\varepsilon\},
\]
whence $C^+_K$ is $\cM_p$-measurable.
\qed

\section*{Aknowledgements}
The authors are very grateful to Mathieu Ribatet and Julien Michel for fruitful discussions and comments on early versions of the manuscript.

\bibliographystyle{plain}
\bibliography{Biblio}

\end{document}